\documentclass[11.5pt]{amsart}

\usepackage{graphicx,amssymb} 
\usepackage{epsfig}

\begin{document}

\pagestyle{myheadings}

\parindent 0mm

\title[A note on the box dimension of degenerate foci]
{A note on the box dimension of degenerate foci}

\author{Domagoj Vlah, Darko \v Zubrini\'c and Vesna \v Zupanovi\'c}

\begin{abstract}
We study polynomial planar systems with singularity of focus type without characteristic directions.
Simple and natural transformation of weak focus has been used to obtain such degenerate focus.
We compute the box dimension of a spiral trajectory, and show connection to cyclicity of the system under a perturbation.
\end{abstract}
\medskip

\date{}

\maketitle

\newtheorem{theorem}{Theorem}
\newtheorem{cor}{Corollary}
\newtheorem{prop}{Proposition}
\newtheorem{lemma}{Lemma}

\font\csc=cmcsc10

\def\esssup{\mathop{\rm ess\,sup}}
\def\essinf{\mathop{\rm ess\,inf}}
\def\wo#1#2#3{W^{#1,#2}_0(#3)}
\def\w#1#2#3{W^{#1,#2}(#3)}
\def\wloc#1#2#3{W_{\scriptstyle loc}^{#1,#2}(#3)}
\def\osc{\mathop{\rm osc}}
\def\var{\mathop{\rm Var}}
\def\supp{\mathop{\rm supp}}
\def\Cap{{\rm Cap}}
\def\norma#1#2{\|#1\|_{#2}}

\def\C{\Gamma}

\let\text=\mbox

\catcode`\@=11
\let\ced=\c
\def\a{\alpha}
\def\b{\beta}
\def\c{\gamma}
\def\d{\delta}
\def\g{\lambda}
\def\o{\omega}
\def\q{\quad}
\def\n{\nabla}
\def\s{\sigma}
\def\div{\mathop{\rm div}}
\def\sing{{\rm Sing}\,}
\def\singg{{\rm Sing}_\ty\,}

\def\A{{\cal A}}
\def\F{{\cal F}}
\def\H{{\cal H}}
\def\W{{\bf W}}
\def\M{{\cal M}}
\def\N{{\cal N}}
\def\S{{\cal S}}

\def\eR{{\bf R}}
\def\eN{{\bf N}}
\def\Ze{{\bf Z}}
\def\Qe{{\bf Q}}
\def\Ce{{\bf C}}

\def\ty{\infty}
\def\e{\varepsilon}
\def\f{\varphi}
\def\:{{\penalty10000\hbox{\kern1mm\rm:\kern1mm}\penalty10000}}
\def\ov#1{\overline{#1}}
\def\D{\Delta}
\def\O{\Omega}
\def\pa{\partial}

\def\st{\subset}
\def\stq{\subseteq}
\def\pd#1#2{\frac{\pa#1}{\pa#2}}
\def\sgn{{\rm sign}\,}
\def\sp#1#2{\langle#1,#2\rangle}

\newcount\br@j
\br@j=0
\def\q{\quad}
\def\gg #1#2{\hat G_{#1}#2(x)}
\def\inty{\int_0^{\ty}}
\def\remark{\smallskip\advance\br@j by1 \noindent{\csc Remark
\the\br@j.}\kern3mm\relax}
\def\example{\smallskip\advance\br@j by1 \noindent{\csc Example
\the\br@j.}\kern3mm\relax}
\def\od#1#2{\frac{d#1}{d#2}}

\def\bg{\begin}
\def\eq{equation}
\def\bgeq{\bg{\eq}}
\def\endeq{\end{\eq}}
\def\bgeqnn{\bg{eqnarray*}}
\def\endeqnn{\end{eqnarray*}}
\def\bgeqn{\bg{eqnarray}}
\def\endeqn{\end{eqnarray}}

\def\bgeqq#1#2{\bgeqn\label{#1} #2\left\{\begin{array}{ll}}
\def\endeqq{\end{array}\right.\endeqn}

\def\abstract{\bgroup\leftskip=2\parindent\rightskip=2\parindent
        \noindent{\bf Abstract.\enspace}}
\def\endabstract{\par\egroup}

\def\udesno#1{\unskip\nobreak\hfil\penalty50\hskip1em\hbox{}
             \nobreak\hfil{#1\unskip\ignorespaces}
                 \parfillskip=\z@ \finalhyphendemerits=\z@\par
                 \parfillskip=0pt plus 1fil}
\catcode`\@=11

\def\cal{\mathcal}
\def\eR{\mathbb{R}}
\def\eN{\mathbb{N}}
\def\Ze{\mathbb{Z}}
\def\Qu{\mathbb{Q}}
\def\Ce{\mathbb{C}}

\section{Introduction}

Analysis of the Poincar\' e map, or the first return map is a standard approach in the study of monodromic singular points and limit cycles. The Poincar\' e map and the normal form  of a weak focus has been studied in \cite{zuzu} and \cite{belg}, from the point of view of fractal geometry. Weak focus is a singular point having pure imaginary eigenvalues of the linear part.
The box dimension of a spiral trajectory near weak focus, and also near a limit cycle has been computed.
Furthermore, the explicit relation between the box dimension and the leading power in the asymptotic expansion of the Poincar\' e map of the weak focus has been obtained.

Here we announce our results which extend the investigation to certain classes of degenerate foci. Nilpotent focus and focus with no linear part are called  degenerate foci. The nilpotent focus is the focus with nilpotent matrix of the linear part. Here we deal with systems having degenerate focus without linear part and without characteristic directions. Characteristic directions can be seen after blowing up of the system.

Here we use results for weak focus from \cite{zuzu} and \cite{belg} and apply them to degenerate focus without characteristic directions. Such degenerate focus has the same  asymptotic of the Poincar\' e map in each direction. In general, the Poincar\' e map of degenerate focus has different asymptotic expansion depending on the direction. That is the reason why the approach to weak focus and to degenerate focus should  be different. Nilpotent focus has two asymptotics, on the characteristic curve and elsewhere. Degenerate focus can have more than one characteristic directions. The asymptotic expansion of the Poincar\' e map near focus has been computed in \cite{medvedeva}. The basic technique, method of blowing-up, shows that these characteristic directions are associated to singular points of the obtained polycycle.

\section{Definitions}

For $A\st\eR^N$  bounded we define the \emph{$\e$-neigh\-bour\-hood} of  $A$ as:
$
A_\e:=\{y\in\eR^N\:d(y,A)<\e\}
$.
By the \emph{lower $s$-dimensional  Minkowski content} of $A$, for $s\ge0$, we mean
$$
\M_*^s(A):=\liminf_{\e\to0}\frac{|A_\e|}{\e^{N-s}},
$$
and analogously for the \emph{upper $s$-dimensional Minkowski content} $\M^{*s}(A)$.
If $\M^{*s}(A)=\M_*^{s}(A)$, we call the common value the \emph{$s$-dimensional Minkowski content of $A$}, and denote it by $\M^s(A)$.
The lower and upper box dimensions of $A$ are
$$
\underline\dim_BA:=\inf\{s\ge0\:\M_*^s(A)=0\}
$$
 and analogously
$\ov\dim_BA:=\inf\{s\ge0\:\M^{*s}(A)=0\}$.
If these two values coincide, we call it simply the box dimension of $A$, and denote it by $\dim_BA$. This will be our situation.
If $0<\M_*^d(A)\le\M^{*d}(A)<\ty$ for some $d$, then we say
 that $A$ is \emph{Minkowski nondegenerate}. In this case obviously $d=\dim_BA$.
In the case when the lower or upper $d$-dimensional Minkowski content of $A$ is equal to $0$ or $\ty$, where $d=\dim_BA$, we say that $A$ is \emph{degenerate}.
If there exists $\M^d(A)$ for some $d$ and $\M^d(A)\in(0,\ty)$, then we say that $A$ is \emph{Minkowski measurable}.

We shall use the following notation.
For any two sequences of positive real numbers $(a_k)$ and $(b_k)$ converging to zero we write $a_k\simeq b_k$ as $k\to\ty$
if there exist positive real numbers $A<B$ such that $a_k/b_k\in[A,B]$ for all $k$.
Also if $f,g:(0,r)\to(0,\ty)$ are two functions converging 
to zero as $s\to0$ and $f(s)/g(s)\in[A,B]$, we write $f(s)\simeq g(s)$
as $s\to0$. We call such sequences and functions comparable.

Spiral trajectory $\Gamma_{\alpha}$ of weak focus  is $\alpha$-\emph{power spiral}, that is spiral $r=f(\f)$ satisfying $f(\f)\simeq\f^{-\a}$ for $0<\alpha\le 1$, see \cite{zuzu}.

\section{$F_{n,n}$ Transformation of weak focus}

We define the map $F_{n,m}$  in order to transform trajectories of weak focus to trajectories of degenerate focus.
Let $F_{n,m}:\eR^2\to\eR^2$ be defined by 
\bgeq\label{F}
F_{m,n}(x,y)=((\sgn x)|x|^{1/m},(\sgn y)|y|^{1/n}).
\endeq
We assume that $m$ and $n$ are positive integers.  Since $F_{1,1}$ is identity, we also assume that 
at least one of these integers is $\ge 2$. It is clear that $F_{m,n}$ is a homeomorphism mapping each quadrant onto itself. Note that it is not diffeomorphism.
We defined the map $F_{m,n}$  in order to transform trajectories of weak focus to trajectories of degenerate focus. A similar idea is used in  e.g. \cite{gt1}, \cite{gt},  \cite{ggg} using the  quasi-homogeneous polar coordinates in the computation of generalized Lyapunov coefficients.

\begin{theorem} 

Assume that $\C$ is a trajectory of the system (\ref{degprvi}), $k,n\ge 1$ 
\bgeq\label{degprvi}
\begin{array}{ccl}
\dot x&=&-y^{2n-1}\pm x^n y^{n-1} (x^{2n}+y^{2n})^k\\
\dot y&=&\phantom{-}x^{2n-1}\pm x^{n-1} y^n (x^{2n}+y^{2n})^k.
\end{array}
\endeq
Then
\bgeq\label{dimprvi}
\dim_B\C=2-\frac{2}{1+2kn},
\endeq
and $\C$ is Minkowski nondegenerate.
\end{theorem}

{\it Idea of the proof.}
Proof is based on application of 
\bgeq
F_{n,n}(x,y)=((\sgn x)|x|^{1/n},(\sgn y)|y|^{1/n})
\endeq
to system 
\bgeq\label{slabi}
\begin{array}{ccl}
\dot x&=&-y\pm x(x^{2}+y^{2})^k\\
\dot y&=&\phantom{-}x\pm y(x^{2}+y^{2})^k
\end{array}
\endeq
with weak focus. For weak focus we apply  Theorem 9 from \cite{zuzu}, as well as the careful study of degenerate focus defined by system (\ref{degprvi}), obtained via (\ref{slabi}).

\remark
A bifurcation parametar $\lambda$ could be added in system (\ref{degprvi}), in order to obtain
\bgeq
\begin{array}{ccl}
\dot x&=&-y^{2n-1}\pm x^n y^{n-1} ((x^{2n}+y^{2n})^k+\lambda)\\
\dot y&=&\phantom{-}x^{2n-1}\pm x^{n-1} y^n((x^{2n}+y^{2n})^k+\lambda).
\end{array}
\endeq
For $\lambda<0$ a limit cycle is born from the degenerate  focus corresponding to parametar $\lambda =0$, and having the box dimension (\ref{dimprvi}).
For $k=1$ and $n=1$ we have standard Hopf bifurcation where a limit cycle has been born from weak focus with $\dim_B\C=4/3$, see \cite{zuzu}.

\medskip

In  \cite{belg} a flow-sector theorem has been proved for weak focus. The theorem says that weak focus flow in sectors near the singular point is lipeomorphically equivalent to the annulus flow. The flow-sector theorem  has been proved for the system with weak focus
\bgeq
\begin{array}{ccl}
\dot x&=&-y+p(x,y)\\
\dot y&=&\phantom{-}x+q(x,y)\nonumber,
\end{array}
\endeq
in which functions $p(x,y)$ and $q(x,y)$ are given $C^1$-functions such that $|p(x,y)|\le C (x^{2}+y^{2})$
and $|q(x,y)|\le C(x^{2}+y^{2})$ for some positive constant $C$ and for $(x,y)$ near the origin.
We generalize flow-sector theorem for the focus of the system (\ref{degprvi}).

For a function $F:U\to V$ with $U,V\st\eR^2$, $V=F(U)$, if  $F$ and $F^{-1}$ are Lipschitzian we say that
$F$ is {\em lipeomorphism}, and that the sets $U$ and $V$ are lipeomorphic. Here the {\em annulus flow} is defined by 
$$
x^{2n}+y^{2n}=const,
$$
for $n\in\mathbb{N}$.

In the following theorem we claim that the focus flow in sectors near the singular point of (\ref{degprvi}) is lipeomorphically equivalent to the annulus flow.

\begin{theorem} \label{flowsectordeg}
Let $U_0\st\eR^2$ be an open sector with
the vertex at the origin, such that its opening angle is in $(0,2\pi)$, and 
the boundary of $U_0$ consists of a part of a trajectory and of intervals on two rays emanating from the origin.
If the diameter of $U_0$ is sufficiently small, then system (\ref{degprvi}) 
restricted to $U_0$ is lipeomorphically equivalent to the system
\bgeq
\begin{array}{ll}
\dot x&=-y^{2n-1}\\
\dot y&=\phantom{-}x^{2n-1}\\\nonumber
\end{array}
\endeq
defined on the sector $V_0$.
\end{theorem}

{\it Idea of the proof.}
The proof is based on the transformation of a weak focus and application of 
 flow sector theorem from \cite{belg}.

\begin{center}
\includegraphics[width=0.5\textwidth]{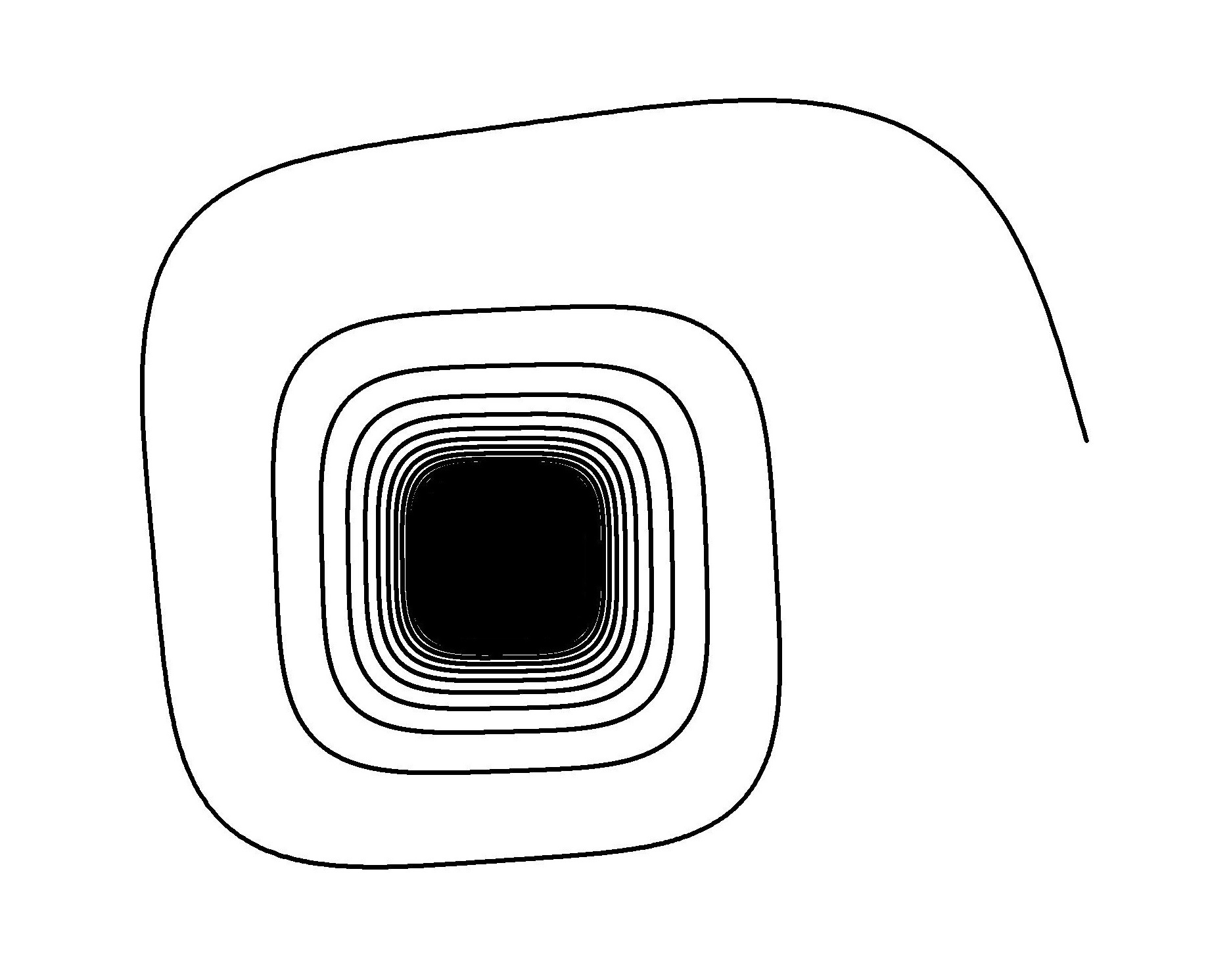}
\end{center}

\remark
It is possible to state previous two theorems for a larger class of systems with degenerate focus without characteristic directions.

\smallskip

Characteristic direction for singularity in the origin of a system
\bgeq\label{general}
\begin{array}{ccl}
\dot x&=& P(x,y)\\
\dot y&=&Q(x,y)\nonumber
\end{array}
\endeq
is linear factor in $\eR [x,y]$ of $yP_d(x,y)-xQ_d(x,y)$, where $P_d(x,y)$ and $Q_d(x,y)$ are homogeneous polynomials of the lowest degree. It is obvious that system (\ref{degprvi}) has no characteristic directions. If there are no characteristic directions, a singular point is either center or focus. The converse is not true.

Now, we consider a different class of degenerate systems.

\begin{prop}
Assume that $\C$ is a trajectory of the system 
\begin{equation}\label{hom}
\begin{array}{ccl}
\dot x&=&-y(x^2+y^2)^{k}+xR_s(x,y)\\
\dot y&=&\phantom{-}x(x^2+y^2)^{k}+yR_s(x,y),
\end{array}
\end{equation}
where $R_s(x,y)$ is a homogeneous polynomial of even degree $s$ and $s> 2k\ge 0$, where $s$ and $k$ are integers.
Then
\bgeq\label{dimex}
d=\dim_B\Gamma=2-\frac{2}{s-2k+1},
\endeq
and $\C$ is Minkowski nondegenerate.
\end{prop}

{\it Idea of the proof.}
Find explicit solution in polar coordinates.

\example
We reveal Example 4 from \cite{ggg} dealing with system (\ref{hom}).  In polar coordinates we obtain 
\begin{equation}
\begin{array}{ccl}
\dot r&=&r^{s+1}R_s(\cos\varphi,\sin\varphi)\\
\dot \varphi&=&r^{2k},
\end{array}
\end{equation}
and the origin is a focus if and only if the following integral is different from zero
\begin{equation}
\int_0^{2\pi}R_s(\cos\varphi ,\sin\varphi)\, d\varphi\ne 0\nonumber.
\end{equation}
We remark that for $s$ odd, the origin is center.
In the focus case a spiral trajectory $\Gamma$ of system (\ref{hom}) is $1/{(s-2k)}$-power spiral with 
\begin{equation}
d=\dim_B\Gamma=2-\frac{2}{s-2k+1}\nonumber.
\end{equation}
Also, box dimension of an orbit generated by the Poincar\' e map, at any transversal through the origin, is equal to $d/2$. 
To system (\ref{hom}) we add analytic perturbation $(\bar P(x,y,\varepsilon), \bar Q(x,y,\varepsilon))=O(\|(x,y)\|)^s,$ satisfying $\bar P(x,y,0)=\bar Q(x,y,0))=0$, for $\varepsilon\in \mathbb{R}$, $0<\|\varepsilon\|\ll 1$. 
Let the corresponding vector field be denoted by $X_{\varepsilon}$. Cyclicity of the origin $p_0$ of the perturbed system will be denoted by $Cycl(X_{\varepsilon},p_0)$. By cyclicity we mean the sharp upper bound for the number of limit cycles which can bifurcate from the origin $p_0$ of the perturbed system.
Theorem 1 from \cite{ggg} says that in this example $Cycl(X_{\varepsilon},p_0)>s/2$. For $k=0$ we obtain weak focus with $Cycl(X_{\varepsilon},p_0)=s/2$.
 According to (\ref{dimex}), we can express cyclicity by using box dimension of the trajectory.

\section{$F_{m,n}$ Transformation of weak focus}

The case of $m\neq n$ is more difficult, and by applying transformation (\ref{F}) to (\ref{slabi}), we obtain the following degenerate system, extending (\ref{degprvi})
\bgeq\label{degdrugi}
\begin{array}{ccl}
\dot x&=&-ny^{2n-1}\pm n x^m y^{n-1} (x^{2m}+y^{2n})^k\\
\dot y&=&\phantom{-}mx^{2m-1}\pm mx^{m-1} y^n (x^{2m}+y^{2n})^k.
\end{array}
\endeq

\begin{theorem}

Assume that $\C$ is a trajectory of the system (\ref{degdrugi}) where $m,n\ge 1$, $m\ge n$. Then 
\bgeq
\dim_B \C\ge 2-\frac{1+\frac{n}{m}}{1+2km}.
\endeq

\end{theorem}

{\it Idea of the proof.}
The proof is based on  the careful study of the degenerate spiral $\C=F_{m,n}(\C_1)$, where $\C_1$ is a spiral trajectory of the system (\ref{slabi}) with weak focus, and $F_{m,n}$  is transformation defined by (\ref{F}).  Also we use suitable Lipschitz transformations.

\remark
We hypothesize that
\bgeq
\dim_B F_{m,n}(\C)= 2-\frac{1+\frac{n}{m}}{1+2km}.
\endeq

\section{Acknowledgments}

This work has been supported in part by Croatian Science Foundation under the project IP-2014-09-2285.


\begin{thebibliography}{12}

\let\small=\rm
\let\normalsize=\rm


\bibitem{gt1} A.\ Gasull, J.\ Torregrosa, A new algorithm for the computation of the Lyapunov constants for some
degenerated critical points, Nonlinear Analysis 47 (2001), 4479--4490.

\bibitem{gt} A.\ Gasull, J.\ Torregrosa, A new approach to the computation of the Lyapunov Constants, Computational and Applied Mathematics, Vol. 20, N. 1-2, (2001), 1-29.


\bibitem{ggg} I.\  Garcia, H.\  Giacomini, M.\  Grau, Journal of Dynamics and Differential Equations, 23, 2,(2011), 251-281


\bibitem{medvedeva} N.\ B. \ Medvedeva,  On the analytic solvability of the problem of distinguishing between center and focus, Proceedings of the Steklov Institute of Mathematics, September 2006, Volume 254, Issue 1, pp 7–93


\bibitem{tricot} C.\ Tricot, {\em Curves and Fractal Dimension}, Springer--Verlag, 1995.


\bibitem{zuzu} D.\ \v Zubrini\'c, V.\ \v Zupanovi\'c, Fractal analysis of spiral trajectories of some
planar vector fields, Bulletin des Sciences Math\'ematiques, 129/6 (2005), 457--485.

\bibitem{belg} D.\ \v Zubrini\'c, V.\ \v Zupanovi\'c, Poincar\'e map in fractal analysis of spiral trajectories of planar vector fields, Bulletin of the Belgian Mathematical Society Simon Stevin,  15 (2008), 5; 947-960 


\end{thebibliography}
\end{document}